\documentclass[12pt]{article}%
\usepackage{amsmath}
\usepackage{amsfonts}
\usepackage{amssymb}
\usepackage{graphicx}
\usepackage{color}%
\providecommand{\U}[1]{\protect\rule{.1in}{.1in}}
\newtheorem{theorem}{Theorem}

\newtheorem{claim}{Claim}

\newtheorem{conjecture}{Conjecture}

\newtheorem{lemma}[theorem]{Lemma}

\newtheorem{observation}[theorem]{Observation}
\newtheorem{problem}{Problem}

\providecommand{\boksie}{\ensuremath{\mathbin{\raisebox{0.3mm}{$\scriptstyle\square$}}}}
\begin{document}

\title{\textbf{Irredundance Trees of Diameter }$\mathbf{3}$}
\author{C.M. Mynhardt\thanks{Supported by the Natural Sciences and Engineering
Research Council of Canada.}\\Department of Mathematics and Statistics\\University of Victoria, Victoria, BC, \textsc{Canada}\\{\small kieka@uvic.ca}
\and A. Roux\thanks{This work is based on research supported by the National
Research Foundation of South Africa (Grant number: 121931)}\\Department of Mathematical Sciences\\Stellenbosch University, Stellenbosch, \textsc{South Africa}\\{\small rianaroux@sun.ac.za}}
\maketitle

\begin{abstract}
A set $D$ of vertices of a graph $G=(V,E)$ is irredundant if each non-isolated
vertex of $G[D]$ has a neighbour in $V-D$ that is not adjacent to any other
vertex in $D$. The upper irredundance number $\operatorname{IR}(G)$ is the
largest cardinality of an irredundant set of $G$; an $\operatorname{IR}%
(G)$-set is an irredundant set of cardinality $\operatorname{IR}(G)$.

The $\operatorname{IR}$-graph of $G$ has the $\operatorname{IR}(G)$-sets as
vertex set, and sets $D$ and $D^{\prime}$ are adjacent if and only if
$D^{\prime}$ can be obtained from $D$ by exchanging a single vertex of $D$ for
an adjacent vertex in $D^{\prime}$. An $\operatorname{IR}$-tree is an
$\operatorname{IR}$-graph that is a tree. We characterize $\operatorname{IR}%
$-trees of diameter $3$ by showing that these graphs are precisely the double
stars $S(2n,2n)$, i.e., trees obtained by joining the central vertices of two
disjoint stars $K_{1,2n}$.

\end{abstract}

\noindent\textbf{Keywords:\hspace{0.1in}}Irredundance; Reconfiguration
problem; $\operatorname{IR}$-graph; $\operatorname{IR}$-tree

\noindent\textbf{AMS Subject Classification Number 2010:\hspace{0.1in}}05C69

\section{Introduction}

\label{Sec_Intro}Reconfiguration problems are concerned with determining
conditions under which a feasible solution to a given problem can be
transformed into another one via a sequence of feasible solutions in such a
way that any two consecutive solutions are adjacent according to a specified
adjacency relation. The solutions form the vertex set of the associated
reconfiguration graph, two vertices being adjacent if one solution can be
obtained from the other in a single step. Typical questions about the
reconfiguration graph concern its structure (connectedness, Hamiltonicity,
diameter), realizability (which graphs can be realized as a specific type of
reconfiguration graph), and algorithmic properties (finding a shortest path
between two solutions).

Domination reconfiguration problems involving dominating sets of different
cardinalities were introduced by Haas and Seyffarth \cite{HS1} and also
studied in, for example, \cite{davood, HS2, haddadan-2015, MRT}. Two
variations of domination reconfiguration problems involving only
minimum-cardinality dominating sets were introduced by Lakshmanan and
Vijayakumar \cite{Laksh} in 2010, and Fricke, Hedetniemi, Hedetniemi, and
Hutson \cite{FHHH11} in 2011, respectively, and also studied in \cite{Bien,
CHH10, Dyck, EdThesis, Laksh, Laura1, Amu, SS08, SMN}. A survey of results
concerning the reconfiguration of colourings and dominating sets in graphs can
be found in~\cite{MN}.

The study of upper irredundance graphs, or $\operatorname{IR}$-graphs, was
mentioned as an open problem in \cite{Laura1} and initiated by the current
authors in \cite{mynhardt19}. There we showed that all disconnected graphs,
but not all connected graphs, are realizable as $\operatorname{IR}$-graphs,
and that the smallest non-complete $\operatorname{IR}$-tree is the double star
$S(2,2)$ (see Figure \ref{Fig_S22}). Here we characterize $\operatorname{IR}%
$-trees of diameter $3$ by showing that these graphs are precisely the
\emph{double stars} $S(2n,2n)$, i.e., trees obtained by joining the central
vertices of two disjoint stars $K_{1,2n}$. We need a number of results from
\cite{mynhardt19}, which we state in Section \ref{Sec_Previous} after
providing some definitions in Section \ref{Sec_Defs}. In Section \ref{Sec_Gn}
we construct a class of graphs $G_{n},\ n\geq1$, and show that
$\operatorname{IR}(G_{n})\cong S(2n,2n)$. We show in Section \ref{Sec_Only}
that these double stars are the only $\operatorname{IR}$-trees of diameter
$3$. We close by mentioning open problems and conjectures in Section
\ref{Sec_Open}.

\section{Definitions}

\label{Sec_Defs}We follow the notation of \cite{CLZ} for general concepts, and
that of \cite{HHS} for domination related concepts not defined here. For a
graph $G=(V,E)$ and vertices $u,v\in V$, we use the notation $u\sim v$
($u\nsim v$, respectively) to denote that $u$ is adjacent (nonadjacent,
respectively) to $v$. For a set $D\subseteq V$ and a vertex $v\in D$, a
$D$-\emph{private neighbour }of $v$ is a vertex $v^{\prime}$ that is dominated
by $v$ (i.e., $v^{\prime}=v$ or $v^{\prime}\sim v$) but by no vertex in
$D-\{v\}$. The set of $D$-private neighbours of $v$ is called the
\emph{private neighbourhood of }$v$ \emph{with respect to }$D$ and denoted by
$\operatorname{PN}(v,D)$.

The concept of irredundance was introduced by Cockayne, Hedetniemi and Miller
\cite{CHM} in 1978. A set $D\subseteq V$ is \emph{irredundant }if
$\operatorname{PN}(v,D)\neq\varnothing$ for each $v\in D$. The \emph{upper
irredundance number }$\operatorname{IR}(G)$ is the largest cardinality of an
irredundant set of $G$. An $\operatorname{IR}$-\emph{set} of $G$, or an
$\operatorname{IR}(G)$-\emph{set}, is an irredundant set of
cardinality$~\operatorname{IR}(G)$. Let $D$ be an irredundant set of $G$. For
$v\in D$, it is possible that $v\in\operatorname{PN}(v,D)$; this happens if
and only if $v$ is isolated in the subgraph $G[D]$ induced by $D$. If
$u\in\operatorname{PN}(v,D)$ and $u\neq v$, then $u\in V-D$; in this case $u$
is an \emph{external }$D$-\emph{private neighbour of }$v$. The set of external
$D$-private neighbours of $v$ is denoted by $\operatorname{EPN}(v,D)$. An
isolated vertex of $G[D]$ may or may not have external $D$-private neighbours,
but if $v$ has positive degree in $G[D]$, then $\operatorname{EPN}%
(v,D)\neq\varnothing$.

\begin{figure}[ptb]%
\centering
\includegraphics[
height=1.4935in,
width=4.1303in
]%
{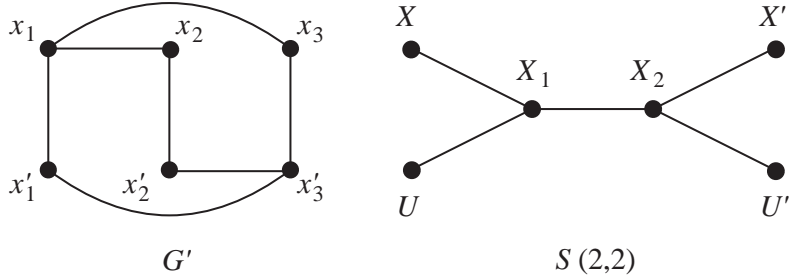}%
\caption{The graph $G^{\prime}=G[\{x_{1},x_{2},x_{3},x_{1}^{\prime}%
,x_{2}^{\prime},x_{3}^{\prime}\}]$ and the double star $S(2,2)$ induced by the
$\operatorname{IR}(G^{\prime})$-sets, as described in Lemma \ref{Lem_Not_Path}%
}%
\label{Fig_S22}%
\end{figure}

As defined in \cite{mynhardt19}, the \emph{$\operatorname{IR}$-graph
}$G(\operatorname{IR})$\emph{ }of \emph{$G$} is the graph whose vertex set
consists of the $\operatorname{IR}(G)$-sets, where sets $D$ and $D^{\prime}$
are adjacent if and only if there exist vertices $u\in D$ and $v\in D^{\prime
}-\{u\}$ such that $uv\in E(G)$ and $D^{\prime}=(D-\{u\})\cup\{v\}$. We
shorten the expression $D^{\prime}=(D-\{u\})\cup\{v\}$ to $D\overset{uv}{\sim
}D^{\prime}$, and also write $D\sim_{H}D^{\prime}$ or $D\overset{uv}{\sim
}_{G(\operatorname{IR})}D^{\prime}$ to show that $D$ and $D^{\prime}$ are
adjacent in $G(\operatorname{IR})$. When $D\overset{uv}{\sim}D^{\prime}$, we
say that $v$ is \emph{swapped into} and $u$ is \emph{swapped out of} the
$\operatorname{IR}(G)$-set, or simply that $u$ and $v$ are \emph{swapped}. To
prove that a given graph $H$ is an $\operatorname{IR}$-graph, one needs to
construct a graph $G$ such that $G(\operatorname{IR})\cong H$. Figure
\ref{Fig_S22} shows a graph $G^{\prime}$ and and its $\operatorname{IR}$-graph
$G^{\prime}(\operatorname{IR})\cong S(2,2)$; the six $\operatorname{IR}%
(G)$-sets are given in Lemma \ref{Lem_Not_Path} below.%

\section{Previous Results}

\label{Sec_Previous}We need the following results from Mynhardt and Roux
\cite{mynhardt19}.\ We weakly partition an $\operatorname{IR}$-set $X$ into
subsets $X^{\operatorname{EPN}}$ and $X^{\mathrm{iso}}$ (one of which may be
empty), where each vertex in $X^{\mathrm{iso}}$ is isolated in $G[X]$ and each
vertex in $X^{\operatorname{EPN}}$ has at least one external private
neighbour. (This partition is not necessarily unique. Isolated vertices of
$G[X]$ with external private neighbours can be allocated arbitrarily to
$X^{\operatorname{EPN}}$ or $X^{\mathrm{iso}}$.) For each $y\in
X^{\operatorname{EPN}}$, let $y^{\prime}\in\operatorname{EPN}(y,X)$ and define
$Y^{\prime}=\{y^{\prime}:y\in X^{\operatorname{EPN}}\}$. Let $X^{\prime
}=(X-X^{\operatorname{EPN}})\cup\{Y^{\prime}\}$; note that $|X|=|X^{\prime}|$.
We call $X^{\prime}$ a \emph{flip-set }of $X$, or to be more precise, the
\emph{flip-set }of $X$ \emph{using }$Y^{\prime}$.

The first result from \cite{mynhardt19} allows us to find more
$\operatorname{IR}$-sets by using external private neighbours in a given
$\operatorname{IR}$-set.

\begin{theorem}
\label{Thm_flip-set}\emph{\cite{mynhardt19}}\hspace{0.1in}If $X$ is an
$\operatorname{IR}(G)$-set, then so is any flip-set of $X$. In particular, if
$X^{\prime}$ is the flip-set\emph{ }of $X$ using $Y^{\prime}$, and $y^{\prime
}\in Y^{\prime}$ belongs to $\operatorname{EPN}(y,X)$, then $y\in
\operatorname{EPN}(y^{\prime},X^{\prime})$.
\end{theorem}

The following lemmas from \cite{mynhardt19} are important tools in the study
of $\operatorname{IR}$-trees.

\begin{lemma}
\label{IRCycle}\emph{\cite{mynhardt19}}\hspace{0.1in}If the $\operatorname{IR}%
$-graph $H$ of $G$ is connected and $X$ is an $\operatorname{IR}(G)$-set such that

\begin{enumerate}
\item[$(i)$] $G[X]$ has exactly one edge, or

\item[$(ii)$] $X$ is independent but at least two vertices have $X$-external
private neighbours,
\end{enumerate}

then $H$ contains an induced $C_{4}$.
\end{lemma}

\begin{lemma}
\label{independent}\emph{\cite{mynhardt19}}\hspace{0.1in}Let $G$ be a graph,
all of whose $\operatorname{IR}$-sets are independent. If the
$\operatorname{IR}$-graph $H$ of $G$ is connected and has order at least
three, then $H$ contains a triangle or an induced $C_{4}$.
\end{lemma}

\begin{lemma}
\label{Lem_Number_IR-sets1}\emph{\cite{mynhardt19}}\hspace{0.1in}If the
$\operatorname{IR}$-graph $H$ of $G$ is connected and $X$ is an
$\operatorname{IR}(G)$-set that contains $k\geq3$ vertices of positive degree
in $G[X]$, or with $X$-external private neighbours, then $\operatorname{diam}%
(H)\geq k$.
\end{lemma}

The last lemma in this section is adapted from Lemma 5.1 of \cite{mynhardt19}.

\begin{lemma}
\label{Lem_Not_Path}Let $H$ be an $\operatorname{IR}$-tree of a graph $G$ such
that $\operatorname{diam}(H)=3$. Suppose $X=\{x_{1},...,x_{r}\}$ is an
$\operatorname{IR}(G)$-set such that exactly three vertices, say $x_{1}%
,x_{2},x_{3}$, have positive degree in $G[X]$. For $i=1,2,3$, let
$x_{i}^{\prime}\in\operatorname{EPN}(x_{i},X)$ and let $X^{\prime}$ be the
flip-set of $X$ using $\{x_{1}^{\prime},x_{2}^{\prime},x_{3}^{\prime}\}$. Note
that $d_{H}(X,X^{\prime})=3$ and assume without loss of generality that
$P:(X=X_{0},X_{1},X_{2},X_{3}=X^{\prime})$, where
\[
X_{1}=\{x_{1}^{\prime},x_{2},...,x_{r}\}\ \ \text{and\ \ }X_{2}=\{x_{1}%
^{\prime},x_{2}^{\prime},x_{3},...,x_{r}\},
\]
is an $X$-$X^{\prime}$ geodesic in $H$. Then

\begin{enumerate}
\item[$(i)$] $G[\{x_{1},x_{2},x_{3},x_{1}^{\prime},x_{2}^{\prime}%
,x_{3}^{\prime}\}]$ is the graph $G^{\prime}$ shown in Figure \ref{Fig_S22};
in particular, the sets $X_{1}$ and $X_{2}$ are independent,

\item[$(ii)$] the set $U=\{x_{1},x_{1}^{\prime},x_{2},x_{4},...,x_{r}\}$ is an
$\operatorname{IR}(G)$-set such that $x_{3}\in\operatorname{EPN}%
(x_{1},U),\ x_{3}^{\prime}\in\operatorname{EPN}(x_{1}^{\prime},U)$ and
$x_{2}^{\prime}\in\operatorname{EPN}(x_{2},U)$,

\item[$(iii)$] denoting the flip-set $\{x_{2}^{\prime},x_{3}^{\prime}%
,x_{3},x_{4},...,x_{r}\}$ of $U$ using $\{x_{3},x_{3}^{\prime},x_{2}^{\prime
}\}$ by $U^{\prime}$, the graph induced by $\{X,X_{1},X_{2},X^{\prime
},U,U^{\prime}\}$ is $S(2,2)$, with edges as shown in Figure \ref{Fig_S22}.
\end{enumerate}
\end{lemma}

We need a few more definitions based on Lemma \ref{Lem_Not_Path}. Observe that
$C=(x_{1},x_{3},x_{3}^{\prime},x_{1}^{\prime},x_{1})$ is an induced $4$-cycle
in $G$, with $x_{1},x_{3}\in X$ and $x_{1}^{\prime},x_{3}^{\prime}\in
X^{\prime}$. Significantly, $x_{1}^{\prime}$ and $x_{3}$ have degree $2$ in
$G[\{x_{1},x_{2},x_{3},x_{1}^{\prime},x_{2}^{\prime},x_{3}^{\prime}\}]$, while
$x_{1}$ and $x_{3}^{\prime}$ have degree $3$. The set $U$ is obtained from $X$
by replacing (not swapping) $x_{3}$ with the vertex $x_{1}^{\prime}$ that is
nonadjacent to it in $C$; imagine a token on $x_{3}$ skipping to
$x_{1}^{\prime}$. We say that $U$ is \emph{the skip-set of} $X$
\emph{replacing }$x_{3}$ \emph{with }$x_{1}^{\prime}$, or simply $U$ is
\emph{a skip-set of }$X$ when the vertices involved are obvious. (Thus $X$ is
the skip-set of $U$ replacing $x_{1}^{\prime}$ with $x_{3}$.) The set
$U^{\prime}$ is not only a flip-set of $U$, but also the skip-set of
$X^{\prime}$ replacing $x_{1}^{\prime}$ with $x_{3}$. We call the set
$\mathcal{X}=\{X,X^{\prime},U,U^{\prime}\}$ of these four related
$\operatorname{IR}(G)$-sets, which are leaves of $H$, a $4$-\emph{cluster}.

\section{Construction of a graph $G_{n}$ such that $G_{n}(\operatorname{IR}%
)\protect\cong S(2n,2n)$}

\label{Sec_Gn}In this section we construct a class of graphs $G_{n}$, where
$n\geq1$, such that $\operatorname{IR}(G_{n})=2n+1$ and $G_{n}%
(\operatorname{IR})$ is the double star $S(2n,2n)$.

\subsection{Construction of $G_{n}$}

\begin{figure}[ptb]%
\centering
\includegraphics[
height=2.4526in,
width=2.8746in
]%
{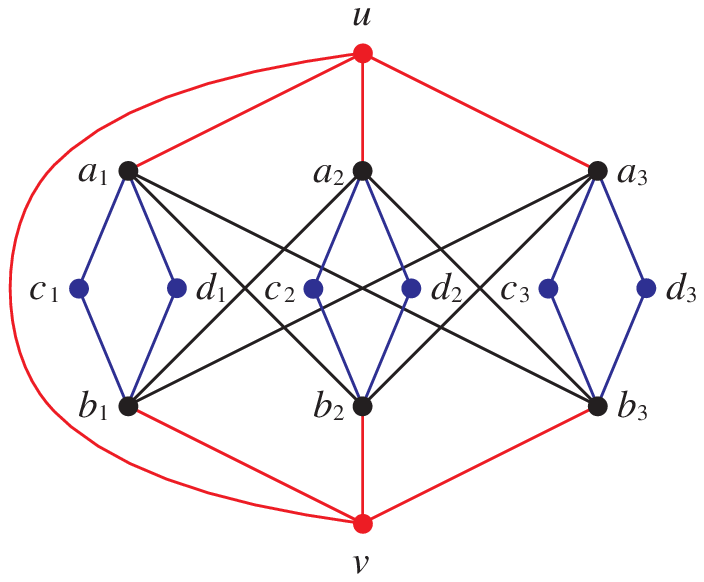}%
\caption{A graph $G_{3}$ such that $G_{3}(\operatorname{IR})\protect\cong
S(6,6)$}%
\label{Fig_Double_star1}%
\end{figure}

See Figure \ref{Fig_Double_star1} for the case $n=3$. The vertex set of
$G_{n}$ consists of five disjoint subsets $\{u,v\},\ A,B,C,D$, where
$A=\{a_{1},...,a_{n}\},\ B=\{b_{1},...,b_{n}\},\ C=\{c_{1},...,c_{n}\}$ and
$D=\{d_{1},...,d_{n}\}$. Join $u$ to each vertex in $A\cup\{v\}$, and join $v$
to each vertex in $B$. Hence, ignoring (for the moment) the edges between $A$
and $B$, the subgraph induced by $A\cup B\cup\{u,v\}$ is isomorphic to the
double star $S(n,n)$. Let $M=\{a_{i}b_{i}:i=1,...,n\}$ and join the vertices
in $A$ to the vertices in $B$ so that the subgraph induced by $A\cup B$ is
isomorphic to the bipartite graph $K_{n,n}-M$ with partite sets $A$ and $B$.
Finally, for each $i=1,...,n$, let $S_{i}=\{a_{i},b_{i},c_{i},d_{i}\}$ and add
edges such that the subgraph induced by each $S_{i}$ is the $4$-cycle
$(a_{i},c_{i},b_{i},d_{i},a_{i})$.%

\subsection{The upper irredundance number of $G_{n}$}

Since $C\cup D\cup\{u\}$ and $C\cup D\cup\{v\}$ are independent sets of
cardinality $2n+1$, $\operatorname{IR}(G_{n})\geq2n+1$. We show that
$\operatorname{IR}(G_{n})=2n+1$. Let $X$ be any irredundant set of $G_{n}$. We
prove that $|X|\leq2n+1$ by proving the following two statements:

\begin{enumerate}
\item[(a)] $|S_{i}\cap X|\leq2$ for each $i$;

\item[(b)] if $\{u,v\}\subseteq X$, then $|X|\leq2n$ for all $n\geq1$.
\end{enumerate}

\noindent\textbf{Proof of (a).\hspace{0.1in}}Note that if, without loss of
generality, $\{a_{i},b_{i},c_{i}\}\subseteq X$ or $\{a_{i},c_{i}%
,d_{i}\}\subseteq X$, then $\operatorname{PN}(c_{i},X)=\varnothing$, a
contradiction.~$\blacklozenge$

\noindent\textbf{Proof of (b).\hspace{0.1in}}Assume that $\{u,v\}\subseteq X$.
Note that $\operatorname{PN}(u,\{u,v\})=A$ and $\operatorname{PN}%
(v,\{u,v\})=B$. Hence, if $n=1$, then $X\cap S_{1}=\varnothing$ and
$X=\{u,v\}$. Assume that $n\geq2$. Suppose $a_{i}\in X$. Then
$\operatorname{PN}(a_{i},X)\subseteq\{c_{i},d_{i}\}$. Moreover,
$\operatorname{PN}(v,X)=\{b_{i}\}$. This implies that $\{b_{i},c_{i}%
,d_{i}\}\cap X=\varnothing$. Since $G_{n}[A\cup B]\cong K_{n,n}-M$, any vertex
$a_{j}$, where\ $j\neq i$, is adjacent to $b_{i}$. Therefore, the private
neighbourhood property of $v$ implies that $A\cap X=\{a_{i}\}$. Similarly, if
$b_{j}\in X$, then $\operatorname{PN}(u,X)=\{a_{j}\},$ $\{a_{j},c_{j}%
,d_{j}\}\cap X=\varnothing$ and $B\cap X=\{b_{j}\}$.

If $\{a_{i},b_{j}\}\subseteq X$ for some $i$ and $j$ (where, as shown above,
$i\neq j$), then possibly $\{c_{k},d_{k}\}\cap X\neq\varnothing$ for some (or
all) $k\notin\{i,j\}$. Necessarily, $|X|\leq4+2(n-2)=2n$. On the other hand,
if $|X\cap(A\cup B)|=1$, say $X\cap(A\cup B)=\{a_{i}\}$, then $X\cap
\{c_{j},d_{j}\}=\varnothing$ for at least one $j\neq i$ to ensure that
$\operatorname{PN}(u,X)\neq\varnothing$. Therefore $|X|\leq3+2(n-2)=2n-1$.

Suppose $X\cap(A\cup B)=\varnothing$. Then $\{c_{i},d_{i}\}\cap X=\varnothing$
for at least one $i$ to ensure that $\operatorname{PN}(u,X)$ and
$\operatorname{PN}(v,X)$ are nonempty. Again we have that $|X|\leq
2+2(n-1)=2n$. This establishes (b).~$\blacklozenge$

\bigskip

Therefore, the cardinality of $X$ is maximized when $|X\cap\{u,v\}|=1$. It
follows that $\operatorname{IR}(G_{n})=2n+1$, as asserted.

\subsection{The $\operatorname{IR}(G_{n})$-sets}

By (a) and (b) and the fact that $\operatorname{IR}(G_{n})=2n+1$, any
$\operatorname{IR}(G_{n})$-set contains exactly one of $u$ and $v$, and
exactly two vertices from each $S_{i}$. We show that the $\operatorname{IR}%
(G_{n})$-sets are precisely the sets

\begin{enumerate}
\item $X=C\cup D\cup\{u\}$

\item $X_{i}=(X-\{c_{i}\})\cup\{a_{i}\}$ for $i=1,...,n$

\item $X_{i}^{\prime}=(X-\{d_{i}\})\cup\{a_{i}\}$ for $i=1,...,n$

\item $Y=C\cup D\cup\{v\}$

\item $Y_{i}=(Y-\{c_{i}\})\cup\{b_{i}\}$ for $i=1,...,n$

\item $Y_{i}^{\prime}=(Y-\{d_{i}\})\cup\{b_{i}\}$ for $i=1,...,n.$
\end{enumerate}

Let $Z$ be any $\operatorname{IR}(G_{n})$-set and assume first that $u\in Z$.
We show that $Z\cap B=\varnothing$. Suppose $b_{i}\in Z$. Since $|Z\cap
S_{i}|=2$, $a_{i}\in Z$ or, without loss of generality, $c_{i}\in Z$.

\begin{itemize}
\item Suppose $\{a_{i},b_{i}\}\subseteq Z$. Note that $N_{G_{n}}%
(u)=A\cup\{v\}$. Since $u\sim a_{i}$, $u$ is not isolated in $G_{n}[Z]$, hence
$u\notin\operatorname{PN}(u,Z)$. Since $a_{i}\in Z$ and $b_{i}$ is adjacent to
$v$ and to each $a_{j},\ j\neq i$, $(A\cup\{v\})\cap\operatorname{PN}%
(u,Z)=\varnothing$. But then $\operatorname{PN}(u,Z)=\varnothing$, which is impossible.

\item Suppose $\{b_{i},c_{i}\}\subseteq Z$. Then $c_{i}$ is not isolated in
$G_{n}[Z]$, hence $c_{i}\notin\operatorname{PN}(c_{i},Z)$. But $b_{i}\in Z$
and $u\sim a_{i}$, hence $\{a_{i},b_{i}\}\cap\operatorname{PN}(c_{i}%
,Z)=\varnothing$. Therefore $\operatorname{PN}(c_{i},Z)=\varnothing$, which is impossible.
\end{itemize}

We conclude that $b_{i}\notin Z$ and therefore $Z\cap B=\varnothing$. We show
next that $|Z\cap A|\leq1$. Suppose to the contrary that $a_{i},a_{j}\in Z$
for $i\neq j$. Since $|Z\cap S_{i}|=2$, we further assume without loss of
generality that $c_{i}\in Z$. Since $a_{i}\sim c_{i}$, $\{a_{i},c_{i}%
\}\cap\operatorname{PN}(c_{i},Z)=\varnothing$, and since $a_{j}\sim b_{i}$,
$b_{i}\notin\operatorname{PN}(c_{i},Z)$. Again we see that $\operatorname{PN}%
(c_{i},Z)=\varnothing$, a contradiction.

We have therefore established that $S_{i}\cap Z=\{c_{i},d_{i}\}$ for all
except possibly one value of $i$. Consider the set $X_{i}=(C-\{c_{i}\})\cup
D\cup\{a_{i},u\}$. For each $j\neq i$, $c_{j}\in\operatorname{PN}(c_{j}%
,X_{i})$ and $d_{j}\in\operatorname{PN}(d_{j},X_{i})$. Since $X_{i}\cap
A=\{a_{i}\}$, no vertex in $X_{i}-\{d_{i}\}$ is adjacent to $b_{i}$, hence
$b_{i}\in\operatorname{PN}(d_{i},X_{i})$. Finally, $v\in\operatorname{PN}%
(u,X_{i})$. Therefore $X_{i}$ is irredundant. Since $|X_{i}|=2n+1$, we have
shown that $X_{i}$ is an $\operatorname{IR}(G_{n})$-set for each $i$.
Similarly, $X_{i}^{\prime}$ is an $\operatorname{IR}(G_{n})$-set for each $i$.
Therefore, the sets $X,\ X_{i}$ and $X_{i}^{\prime}$ are precisely the
$\operatorname{IR}(G_{n})$-sets containing $u$. By symmetry, $Y,\ Y_{i}$ and
$Y_{i}^{\prime}$ are precisely the $\operatorname{IR}(G_{n})$-sets
containing~$v$. Note that $\{X_{i},X_{i}^{\prime},Y_{i},Y_{i}^{\prime}\}$ is a
$4$-cluster for each $i$.%
\begin{figure}[ptb]%
\centering
\includegraphics[
height=1.7521in,
width=2.9308in
]%
{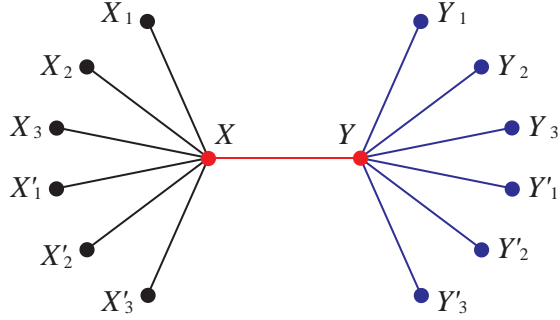}%
\caption{The graph $G_{3}(\operatorname{IR})$, where $G_{3}$ is the graph in
Figure \ref{Fig_Double_star1}, $X=\{u,c_{1},c_{2},c_{3},d_{1},d_{2}%
,d_{3}\},\ Y=\{v,c_{1},c_{2},c_{3},d_{1},d_{2},d_{3}\},\ X_{i}=(X-\{c_{i}%
\})\cup\{a_{i}\},\ X_{i}^{\prime}=(X-\{d_{i}\})\cup\{a_{i}\},\ Y_{i}%
=(Y-\{c_{i}\})\cup\{b_{i}\}$, and $Y_{i}^{\prime}=(Y-\{d_{i}\})\cup\{b_{i}\}$}%
\label{Fig_Double_star2}%
\end{figure}

\subsection{The graph $G_{n}(\operatorname{IR})$}

Denote the $\operatorname{IR}$-graph $G_{n}(\operatorname{IR})$ by $H$. Figure
\ref{Fig_Double_star2} shows $H$ for the case $n=3$. Since $u\sim_{G_{n}}v$,
we see that $X\overset{uv}{\sim}_{H}Y$. Since $c_{i}\sim_{G_{n}}a_{i}$ and
$d_{i}\sim_{G_{n}}a_{i}$, it follows that $X_{i}\overset{c_{i}a_{i}}{\sim}%
_{H}X$ and $X_{i}^{\prime}\overset{d_{i}a_{i}}{\sim}_{H}X$ for each $i$.
Similarly, $Y_{i}\overset{c_{i}b_{i}}{\sim}_{H}Y$ and $Y_{i}^{\prime
}\overset{d_{i}b_{i}}{\sim}_{H}Y$ for each $i$. This shows that the double
star $S(2n,2n)$ is a spanning subgraph of $H$. It remains to show that the
vertices $X_{i},\ X_{i}^{\prime},\ Y_{i}$ and $Y_{i}^{\prime}$ are leaves
of~$H$.

\begin{itemize}
\item When $i\neq j$, note that $\{a_{i},c_{j}\}\subseteq X_{i}-X_{j}$.
Therefore $|X_{i}-X_{j}|\geq2$, hence $X_{i}\nsim_{H}X_{j}$. Similarly,
$X_{i}^{\prime}\nsim_{H}X_{j}^{\prime},\ Y_{i}\nsim_{H}Y_{j}$ and
$Y_{i}^{\prime}\nsim_{H}Y_{j}^{\prime}$.

\item For each $i$, since $c_{i}\nsim d_{i}$, it follows that $X_{i}\nsim
_{H}X_{i}^{\prime}$. Similarly, $Y_{i}\nsim_{H}Y_{i}^{\prime}$.

\item Consider $X_{i}$ and $X_{j}^{\prime}$ for $i\neq j$. Since
$\{a_{i},d_{j}\}\subseteq X_{i}-X_{j}^{\prime}$, we know that $|X_{i}%
-X_{j}^{\prime}|\geq2$ and thus $X_{i}\nsim_{H}X_{j}^{\prime}$. Similarly,
$Y_{i}\nsim_{H}Y_{j}^{\prime}$.
\end{itemize}

It follows that the sets
\[
\mathcal{X}=\{X_{i}:i=1,...,n\}\cup\{X_{i}^{\prime}:i=1,...,n\}
\]
and%
\[
\mathcal{Y}=\{Y_{i}:i=1,...,n\}\cup\{Y_{i}^{\prime}:i=1,...,n\}
\]
are independent in~$G_{n}$.

\begin{itemize}
\item For any $U\in\mathcal{X}$ and any $W\in\mathcal{Y}$, $\{a_{i}%
,u\}\subseteq U-W$ and $\{a_{i},u\}\subseteq U-Y$. Hence no vertex in
$\mathcal{X}$ is adjacent, in $G_{n}(\operatorname{IR})$, to $Y$ or to any
$W\in\mathcal{Y}$.

\item Finally, no vertex $W\in\mathcal{Y}$ is adjacent to $X$, because
$\{b_{i},v\}\subseteq W-X$.
\end{itemize}

This completes the proof that $H=G_{n}(\operatorname{IR})\cong S(2n,2n)$.

\section{Stars $S(2k,2k)$ are the only $\operatorname{IR}$-trees of diameter
$3$}

\label{Sec_Only}Having shown in Section \ref{Sec_Gn} that $S(2n,2n)$ is an
$\operatorname{IR}$-tree for each $n\geq1$, we now show that they are the only
$\operatorname{IR}$-trees of diameter $3$. We begin by stating a simple
observation about trees of diameter $3$ for referencing.

\begin{observation}
\label{Diam3_tree}Let $T$ be a tree with diameter $3$ and diametrical path
$(v_{0},v_{1},v_{2},v_{3})$. For any vertex $u$ of $T$, $d(u,v_{1})\leq2$ and
$d(u,v_{2})\leq2$. Moreover, if $d(u,v_{1})=2$, then $v_{2}$ lies on the
$v_{1}$-$u$ geodesic and $u\sim v_{2}$.
\end{observation}

We are now ready to prove our main result. We use the graphs in Figure
\ref{Fig_S22} throughout the proof to affirm adjacencies and nonadjacencies
between their vertices. To prove the theorem we show that if a tree $T$ with
$\operatorname{diam}(T)=3$ is an $\operatorname{IR}$-tree, then the leaves of
$T$ occur in disjoint $4$-clusters. We state and prove several claims within
the main proof; the end of the proof of a claim is signalled by a diamond
$\blacklozenge$. We consider three cases in the proof of Claim \ref{Cl4}; the
end of the proof of each case is indicated by an open diamond $\lozenge$.

\begin{theorem}
\label{Thm_Main}The stars $S(2k,2k)$, $k\geq1$, are the only
$\operatorname{IR}$-trees of diameter $3$.
\end{theorem}

\noindent\textbf{Proof}.\hspace{0.1in}Suppose $T$ with $\operatorname{diam}%
(T)=3$ is an $\operatorname{IR}$-tree of a graph $G$. By Lemma~\ref{IRCycle},
all the $\operatorname{IR}$-sets of $G$ are either independent or induce a
graph that has exactly three vertices of positive degree. If $G$ has only
independent $\operatorname{IR}$-sets, then by Lemma~\ref{independent} the
$\operatorname{IR}$-graph of $G$ contains a cycle. Therefore $G$ has an
$\operatorname{IR}$-set $X=\{x_{1},...,x_{r}\}$ with exactly three vertices,
say $x_{1},x_{2},x_{3}$, of positive degree in $G[X]$. Let $x_{i}^{\prime}%
\in\operatorname{EPN}(x_{i},X)$ for $i=1,2,3$ and let $X^{\prime}$ be the
flip-set of $X$ using $\{x_{1}^{\prime},x_{2}^{\prime},x_{3}^{\prime}\}$.
Since $|X-X^{\prime}|=3$, we know that $d_{T}(X,X^{\prime})\geq3$, and since
$\operatorname{diam}(T)=3$ by assumption, $d_{T}(X,X^{\prime})=3$. Let
$P:(X=X_{0},X_{1},X_{2},X_{3}=X^{\prime})$ be an $X$-$X^{\prime}$ geodesic in
$T$. Assume, as in the statement of Lemma \ref{Lem_Not_Path}, that
\[
X_{1}=\{x_{1}^{\prime},x_{2},x_{3},\dots,x_{r}\},X_{2}=\{x_{1}^{\prime}%
,x_{2}^{\prime},x_{3},\dots,x_{r}\}\ \text{and\ }X_{3}=X^{\prime}%
=\{x_{1}^{\prime},x_{2}^{\prime},x_{3}^{\prime},\dots,x_{r}\}.
\]
Let
\begin{align*}
Z  &  =\{x_{1},x_{2},x_{3},x_{1}^{\prime},x_{2}^{\prime},x_{3}^{\prime}\},\\
U  &  =\{x_{1},x_{1}^{\prime},x_{2},x_{4},...,x_{r}\},\text{\ }U^{\prime
}=\{x_{2}^{\prime},x_{3}^{\prime},x_{3},x_{4},...,x_{r}\},\text{\ and}\\
\mathcal{X}  &  =\{X,X^{\prime},U,U^{\prime}\}.
\end{align*}
By Lemma \ref{Lem_Not_Path}, $G[Z]$ is the graph $G^{\prime}$ shown in Figure
\ref{Fig_S22}, the sets $X_{1}$ and $X_{2}$ are independent,
\begin{equation}
x_{3}\in\operatorname{EPN}(x_{1},U),\ x_{3}^{\prime}\in\operatorname{EPN}%
(x_{1}^{\prime},U)\text{\ and\ }x_{2}^{\prime}\in\operatorname{EPN}(x_{2},U),
\label{eq_epns}%
\end{equation}
$U^{\prime}$ is the flip-set of $U$ using $\{x_{3},x_{3}^{\prime}%
,x_{2}^{\prime}\}$, and $\mathcal{X}$ is a $4$-cluster. Also by Lemma
\ref{Lem_Not_Path},%
\[
T[\{X,X_{1},X_{2},X^{\prime},U,U^{\prime}\}]\cong S(2,2),
\]
with edges shown in Figure \ref{Fig_S22}. For referencing we state and prove
the following claim.

\begin{claim}
\label{Cl1}The vertex $x_{2}\in\operatorname{EPN}(x_{2}^{\prime},X_{2})$ and
$x_{2}^{\prime}$ is the only vertex with an $X_{2}$-external private neighbour.
\end{claim}

\noindent\textbf{Proof of Claim \ref{Cl1}.\hspace{0.1in}}Since $X_{1}$ is
independent, $x_{1}^{\prime}\nsim x_{2}\nsim x_{3}$. Since $x_{i}$ is isolated
in $G[X]$ for $i=4,...,r$, $x_{2}\nsim x_{i}$ for $i=4,...,r$. Finally, since
$x_{2}^{\prime}\sim x_{2}$, it follows that $x_{2}\in\operatorname{EPN}%
(x_{2}^{\prime},X_{2})$. Since $X_{2}$ is independent, Lemma~\ref{IRCycle}
implies that that no other vertex in $X_{2}$ has an $X_{2}$-external private
neighbour.~$\blacklozenge$\medskip

By assumption, $T$ is a tree, $\operatorname{diam}(T)=3$ and $P$ is a
diametrical path of $T$. Hence any vertex of $T$ not on $P$ is a leaf adjacent
to $X_{1}$ or $X_{2}$. Let $A^{\prime}\notin\{X_{1},X^{\prime},U^{\prime}\}$
be any $\operatorname{IR}(G)$-set such that $X_{2}\sim_{T}A^{\prime}$; say
$X_{2}\overset{yy^{\prime}}{\sim}_{T}A^{\prime}$ for vertices $y\in X_{2}$ and
$y^{\prime}\in V(G)-X_{2}$ adjacent to $y$. We aim to show that $A^{\prime}$
belongs to a $4$-cluster disjoint from $\mathcal{X}$. Since $X_{2}$ is
independent and $y\sim y^{\prime}$, we know that
\begin{equation}
y\in\operatorname{EPN}(y^{\prime},A^{\prime}). \label{eq_y'}%
\end{equation}
We further investigate the neighbourhoods of $y$ and $y^{\prime}$ in Claims
\ref{Cl2} -- \ref{Cl5}.

\begin{claim}
\label{Cl2}The vertex $y^{\prime}$ is adjacent to exactly two vertices of
$A^{\prime}$.
\end{claim}

\noindent\textbf{Proof of Claim \ref{Cl2}.\hspace{0.1in}}First suppose that
$y^{\prime}$ is an $X_{2}$-private neighbour of $y$. Then, by Claim \ref{Cl1},
$y=x_{2}^{\prime}$ (but $y^{\prime}\neq x_{2}$, otherwise $A^{\prime}=X_{1}$)
and $A^{\prime}=\{x_{1}^{\prime},y^{\prime},x_{3},\dots,x_{r}\}$. Furthermore,
if $y^{\prime}\nsim x_{2}$, then $Q_{1}=\{x_{1}^{\prime},y^{\prime}%
,x_{2},x_{3},\dots,x_{r}\}$ is an independent set, hence an irredundant set,
with $X_{1}\subsetneqq Q_{1}$. But then $|Q_{1}|>\operatorname{IR}G)$, which
is impossible. Hence $y^{\prime}\sim x_{2}$. However, now $(X_{1}%
,X_{2},A^{\prime},X_{1})$ is a cycle in the tree $T$. This contradiction
implies that $y^{\prime}$ is not an $X_{2}$-external private neighbour of $y$
and therefore $y^{\prime}$ is not isolated in $G[A^{\prime}]$. In fact, since
$X_{2}$ is independent, it follows from Lemma~\ref{IRCycle} that $y^{\prime}$
is adjacent to at least two vertices in $A^{\prime}=(X_{2}-\{y\})\cup
\{y^{\prime}\}$. If $y^{\prime}$ is adjacent to more than two vertices in
$A^{\prime}$, then the flip-set of $A^{\prime}$ through $N[y^{\prime}]\cap
A^{\prime}$ is at distance at least $3$ from $X_{2}$, contradicting
Observation \ref{Diam3_tree}. We conclude that $y^{\prime}$ is adjacent to
exactly two vertices of $A^{\prime}$.~$\blacklozenge$\medskip

Since $A^{\prime}\notin\{X_{1},X^{\prime},U^{\prime}\}$ is an arbitrary
$\operatorname{IR}(G)$-set adjacent to $X_{2}$ (and since $X^{\prime}$ and
$U^{\prime}$ are not independent), it follows that no $\operatorname{IR}%
(G)$-set adjacent to $X_{2}$, other than $X_{1}$, is independent. By symmetry,
the same statement is true for $\operatorname{IR}(G)$-sets, other than $X_{2}%
$, adjacent to $X_{1}$. Since $\operatorname{diam}(T)=3$ and $P$ is an
$X$-$X^{\prime}$ geodesic, all $\operatorname{IR}(G)$-sets are adjacent to
either $X_{1}$ or $X_{2}$. We therefore deduce that
\begin{equation}
X_{1}\text{\ and\ }X_{2}\text{\ are\ the\ only\ independent\ }%
\operatorname{IR}(G)\text{-sets.} \label{eq_indep}%
\end{equation}

Let $u,v\in N(y^{\prime})\cap A^{\prime}$ with $A^{\prime}$-external private
neighbours $u^{\prime}$ and $v^{\prime}$, respectively. By (\ref{eq_y'}),
$y\in\operatorname{EPN}(y^{\prime},A^{\prime})$. Denote the flip-set of
$A^{\prime}$ using $\{y,u^{\prime},v^{\prime}\}$ by $A$.

\begin{claim}
\label{Cl3}The vertex $y^{\prime}$ belongs to $V(G)-(X\cup Z)$\emph{.}
\end{claim}

\noindent\textbf{Proof of Claim \ref{Cl3}.\hspace{0.1in}}Suppose the vertex
$y\in X_{2}$ is swapped out for one of the vertices $x_{1},x_{2},x_{3}%
^{\prime}$. Recall that, by (\ref{eq_epns}), $x_{i}^{\prime}\in
\operatorname{EPN}(x_{i},X)$ for $i\in\{1,2,3\}$. Since the only neighbours of
$x_{1}$ in $X$ are $x_{2}$ and $x_{3}$, the only neighbours of $x_{1}$ in
$X_{2}$ are $x_{1}^{\prime}$ and $x_{3}$. Since the only neighbour of $x_{2}$
in $X$ is $x_{1}$, the only neighbour of $x_{2}$ in $X_{2}$ is $x_{2}^{\prime
}$. The only neighbours of $x_{3}^{\prime}$ in $X_{2}$ are $x_{3}$ and, by
Lemma \ref{Lem_Not_Path}, $x_{1}^{\prime}$ and $x_{2}^{\prime}$ (see Figure
\ref{Fig_S22}). Hence $y\in\{x_{1}^{\prime},x_{2}^{\prime},x_{3}\}$. But by
Lemma \ref{Lem_Not_Path}, the sets $X_{1},X_{2}$ and the $4$-cluster
$\{X,X^{\prime},U,U^{\prime}\}$ are the only $\operatorname{IR}(G)$-sets that
can be obtained in this way. We conclude that $y\in X_{2}$ is swapped for
$y^{\prime}\in V(G)-X-Z$.~$\blacklozenge$

\begin{claim}
\label{Cl4}The vertex $y=x_{j}$ for some $j\geq4$.
\end{claim}

\noindent\textbf{Proof of Claim \ref{Cl4}.\hspace{0.1in}}We showed in the
proof of Claim \ref{Cl3} that no vertex in $X_{2}$ is swapped out for one of
the vertices $x_{1},x_{2},x_{3}^{\prime}$. Now suppose that $y\in
\{x_{1}^{\prime},x_{2}^{\prime},x_{3}\}$ and $y$ is swapped for $y^{\prime
}\notin Z$. We consider three cases, depending on $y$. In each case we examine
the possibilities for the external $A^{\prime}$-private neighbours $u^{\prime
}$ and $v^{\prime}$ of the neighbours $u$ and $v$, respectively, of
$y^{\prime}$ in $A^{\prime}$. We show that each choice where $\{u^{\prime
},v^{\prime}\}\cap Z\neq\varnothing$ leads to a contradiction. We are then
left with a flip-set $A$ of $A^{\prime}$ containing the vertices $u^{\prime
},v^{\prime}\neq x_{i},x_{i}^{\prime}$ for $i=1,2,3$, which implies that
$d_{T}(X_{i},A)\geq2$ for $i=1,2$, contrary to Observation \ref{Diam3_tree}.
\smallskip

\noindent\textbf{Case 1:\hspace{0.1in}}$y=x_{2}^{\prime}$, that is,
$A^{\prime}=\{x_{1}^{\prime},y^{\prime},x_{3},\dots,x_{r}\}$. We show that
$u^{\prime},v^{\prime}\neq x_{i},x_{i}^{\prime}$ for $i=1,2,3$. Since
$u^{\prime},v^{\prime}$ are external $A^{\prime}$-private neighbours and
$x_{1}^{\prime},x_{3}\in A^{\prime}$, we know that $u^{\prime},v^{\prime}\neq
x_{1}^{\prime},x_{3}$. Since $y^{\prime}\sim x_{2}^{\prime}=y$, no neighbour
of $y^{\prime}$ has $x_{2}^{\prime}$ as $A^{\prime}$-private neighbour, hence
$u^{\prime},v^{\prime}\neq x_{2}^{\prime}$. Also, since $x_{2}\in
\operatorname{EPN}(x_{2}^{\prime},X_{2})$, we either have $N[x_{2}]\cap
A^{\prime}=\varnothing$ or $N[x_{2}]\cap A^{\prime}=\{y^{\prime}\}$; in either
case, $u^{\prime},v^{\prime}\neq x_{2}$. Because $x_{1}^{\prime}\sim
x_{3}^{\prime}\sim x_{3}$, we know that $u^{\prime},v^{\prime}\neq
x_{3}^{\prime}$. Finally, since $x_{1}^{\prime}\sim x_{1}\sim x_{3}$, we have
$u^{\prime},v^{\prime}\neq x_{1}$. Therefore the flip-set $A$ of $A^{\prime}$
contains the vertices $u^{\prime},v^{\prime}\notin Z$, which implies that
$\min\{d_{T}(A,X_{1}),d_{T}(A,X_{2})\}\geq2$ and leads to a contradiction to
Observation \ref{Diam3_tree}.$~\lozenge$\smallskip

\begin{figure}[ptb]%
\centering
\includegraphics[
height=1.9026in,
width=4.1303in
]%
{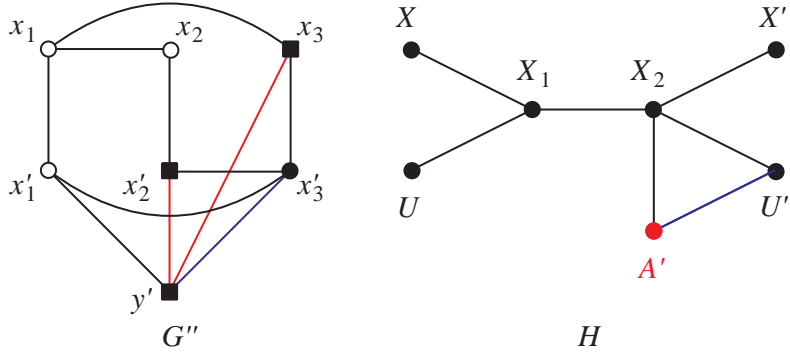}%
\caption{The graph $G^{\prime\prime}=G[Z\cup\{y^{\prime}\}]$ and the graph $H$
formed by its $\operatorname{IR}$-sets, as described in Case 2 of the proof of
Claim \ref{Cl4}}%
\label{Fig_S22_plus_y}%
\end{figure}

\noindent\textbf{Case 2:\hspace{0.1in}}$y=x_{1}^{\prime}$, that is,
$A^{\prime}=\{y^{\prime},x_{2}^{\prime},x_{3},...,x_{r}\}$. We consider the
possibilities for the neighbours $u,v$ of $y^{\prime}$ in $A^{\prime}$ and
their respective external $A^{\prime}$-private neighbours $u^{\prime
},v^{\prime}$. If $y^{\prime}\sim x_{3}^{\prime}$ (see the black and blue
edges in Figure \ref{Fig_S22_plus_y}), then $A^{\prime}\overset{y^{\prime
}x_{3}^{\prime}}{\sim}_{T}U^{\prime}$, so that the graph $H$ in Figure
\ref{Fig_S22_plus_y} is a subgraph of the tree $T$, which is impossible.%
\begin{equation}
\text{Hence\ we\ may\ assume\ that\ }y^{\prime}\nsim x_{3}^{\prime}.
\label{eq_y'x'3}%
\end{equation}
Clearly, $u^{\prime},v^{\prime}\neq x_{2}^{\prime},x_{3}\in A^{\prime}$
because $u^{\prime}$ and $v^{\prime}$ are external $A^{\prime}$-private neighbours.

We next show that $y^{\prime}\sim x_{2}^{\prime}$. Suppose, to the contrary,
that $y^{\prime}\nsim x_{2}^{\prime}$, that is, $u,v\in\{x_{3},...,x_{r}\}$.
We may then assume that $v\neq x_{3}$. By the private neighbour properties of
$X$, this implies that $v^{\prime}\neq x_{i},x_{i}^{\prime}$ for $i=1,2,3$,
since none of these vertices is adjacent to a vertex in $\{x_{4},...,x_{r}\}$.
If we also have that $u\neq x_{3}$, the same holds for $u^{\prime}$, i.e.,
$u^{\prime},v^{\prime}\notin Z$, and we are done, hence assume $u=x_{3}$.
Since $x_{2}^{\prime}\sim x_{3}^{\prime}$, we have that $x_{3}^{\prime}%
\notin\operatorname{EPN}(x_{3},A^{\prime})\cup\operatorname{EPN}(y^{\prime
},A^{\prime})$. Therefore $u^{\prime}\neq x_{3}^{\prime}$. Since $x_{2}\sim
x_{2}^{\prime}$, we know that $u^{\prime}\neq x_{2}$. Suppose $u^{\prime
}=x_{1}$. Then $A=\{x_{1},x_{1}^{\prime},x_{2}^{\prime},v^{\prime}%
,x_{4},...,x_{r}\}-\{v\}$. However, now $x_{1},v^{\prime}\in A-(X_{1}\cup
X_{2})$, so that $\min\{d_{T}(A,X_{1}),d_{T}(A,X_{2})\}\geq2$, contrary to
Observation \ref{Diam3_tree}. We obtain a similar contradiction if $u^{\prime
}\notin Z$.
\begin{equation}
\text{Therefore\ we\ may\ assume\ that}\ y^{\prime}\sim x_{2}^{\prime
}\text{;\ say\ }x_{2}^{\prime}=v. \label{eq_y'x'2}%
\end{equation}
Because $x_{2}^{\prime}\sim x_{3}^{\prime}\sim x_{3}$, we know that
$x_{3}^{\prime}\notin\{u^{\prime},v^{\prime}\}$. We still need to consider
$x_{1}$ and $x_{2}$. Observe that
\begin{equation}
\text{if\ }x_{2}\in\{u^{\prime},v^{\prime}\}\text{,\ then\ }x_{2}=v^{\prime}
\label{eq_x2}%
\end{equation}
because $x_{2}\sim x_{2}^{\prime}=v$ by (\ref{eq_y'x'2}). Since $x_{2}%
^{\prime}\in\operatorname{EPN}(x_{2},X)$, $x_{1}\nsim x_{2}^{\prime}=v$ and so
$x_{1}\neq v^{\prime}$. The investigation of whether $x_{1}=u^{\prime}$
depends on whether $u=x_{3}$ or $u=x_{j}$ for $j\geq4$.\smallskip%

\noindent\textbf{Case 2.1:\hspace{0.1in}}$u=x_{3}$. Consider $A=\{y,u^{\prime
},v^{\prime},x_{4},..,x_{r}\}$.

\noindent Suppose that $u^{\prime}=x_{1}$ and $v^{\prime}=x_{2}$. Then
\[
x_{1}^{\prime}\in\operatorname{EPN}(y^{\prime},A^{\prime}),\ \ x_{2}%
\in\operatorname{EPN}(x_{2}^{\prime},A^{\prime})\ \ \text{and}\ \ x_{1}%
\in\operatorname{EPN}(x_{3},A^{\prime}).
\]
At this point the graph $G[Z\cup\{y^{\prime}\}]$ is the graph formed by the
black and red edges of the graph $G^{\prime\prime}$ in Figure
\ref{Fig_S22_plus_y}, where the vertices in $A^{\prime}-\{x_{4},...,x_{r}\}$
and their private neighbours are shown by solid squares and open circles,
respectively. By Claim \ref{Cl2}, (\ref{eq_y'x'2}) and the assumption of Case
2.1, $y^{\prime}$ is adjacent to exactly two vertices in $A^{\prime}$, namely
$x_{2}^{\prime}$ and $x_{3}$; we deduce that $y^{\prime}\in\operatorname{EPN}%
(x_{3},X)$. Since, by (\ref{eq_y'x'3}), $y^{\prime}\nsim x_{3}^{\prime}$, it
now follows that the set $Q_{2}=(X-\{x_{1},x_{3}\})\cup\{y^{\prime}%
,x_{3}^{\prime}\}=\{y^{\prime},x_{2},x_{3}^{\prime},x_{4},...,x_{r}\}$ is
independent and therefore an $\operatorname{IR}(G)$-set, contradicting
(\ref{eq_indep}).

\begin{figure}[ptb]%
\centering
\includegraphics[
height=1.9026in,
width=4.2843in
]%
{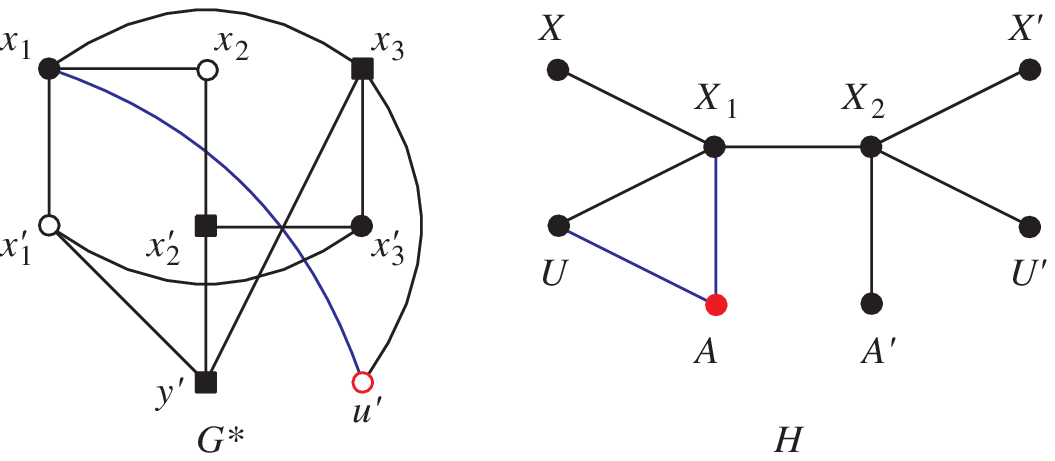}%
\caption{The graph $G^{\ast}=G[Z\cup\{u^{\prime},y^{\prime}\}]$ and the graph
$H$ formed by its $\operatorname{IR}$-sets, as described in Case 2.1 of the
proof of Theorem \ref{Thm_Main}.}%
\label{Fig_S22_plus_u}%
\end{figure}

Suppose $v^{\prime}=x_{2}$ but $x_{1}\notin\operatorname{EPN}(x_{3},A^{\prime
})$, i.e., $x_{1}\neq u^{\prime}$. Then $u^{\prime}\in V(G)-(X\cup Z)$. (See
the black edges of $G^{\ast}$ in Figure \ref{Fig_S22_plus_u}.) Consider the
set $Q_{3}=\{x_{1},x_{2}^{\prime},x_{3},...,x_{r}\}$ and note that $x_{1}%
x_{3}$ is the only edge of $G[Q_{3}]$. Moreover, $x_{1}^{\prime}%
\in\operatorname{EPN}(x_{1},Q_{3})$ and (with only the edges in black)
$u^{\prime}\in\operatorname{EPN}(x_{3},Q_{3})$. But by Lemma~\ref{IRCycle},
$Q_{3}$ is not an $\operatorname{IR}(G)$-set. The only possibility is that
$u^{\prime}$ is adjacent to some other vertex of $Q_{3}$. Since $x_{2}%
^{\prime}\in A^{\prime}$ and $u^{\prime}\in\operatorname{EPN}(x_{3},A^{\prime
})$, it follows that $u^{\prime}\nsim x_{2}^{\prime}$, and we conclude that
$u^{\prime}\sim x_{1}$. (See the black and blue edges of $G^{\ast}$ in Figure
\ref{Fig_S22_plus_u}.) But now $A=\{x_{1}^{\prime},x_{2},u^{\prime}%
,x_{4},...,x_{r}\}$ belongs to the triangle $A\overset{u^{\prime}x_{3}}{\sim
}_{T}X_{1}\overset{x_{3}x_{1}}{\sim}_{T}U\overset{x_{1}u^{\prime}}{\sim}_{T}A$
in $T$, which is impossible.

Finally, suppose $v^{\prime}\neq x_{2}$. If $u^{\prime}=x_{1}$, then
$A=\{x_{1},x_{1}^{\prime},v^{\prime},x_{4},...,x_{r}\}$, where $v^{\prime}\neq
x_{i},x_{i}^{\prime}$ for $i=1,2,3$. Hence $\{x_{1},v^{\prime}\}\subseteq
A-(X_{1}\cup X_{2})$ and we get a contradiction to Observation
\ref{Diam3_tree}. We obtain a similar contradiction if $u^{\prime}\neq x_{1}$,
because then $u^{\prime},v^{\prime}\neq x_{i},x_{i}^{\prime}$ for $i=1,2,3$
and $\{u^{\prime},v^{\prime}\}\subseteq A-(X_{1}\cup X_{2})$. This concludes
the proof of Case 2.1 ($u=x_{3}$).\smallskip%

To recapitulate, we may assume, by (\ref{eq_y'x'3}) and (\ref{eq_y'x'2}), that
$y^{\prime}\nsim x_{3}^{\prime}$ and $y^{\prime}\sim x_{2}^{\prime}=v$. We
have shown that $\{u^{\prime},v^{\prime}\}\cap\{x_{2}^{\prime},x_{3}%
,x_{3}^{\prime}\}=\varnothing$, and if $x_{2}\in\{u^{\prime},v^{\prime}\}$,
then $x_{2}=v^{\prime}$ -- see (\ref{eq_x2}).\smallskip

\noindent\textbf{Case 2.2:\hspace{0.1in}}$u\in\{x_{4},...,x_{r}\}$; say
$u=x_{4}$. By the properties of $X$, $u^{\prime}\notin X\cup Z$; say
$u^{\prime}=x_{4}^{\prime}$. Moreover, since $x_{1}\nsim x_{2}^{\prime}$, we
have that $x_{1}\neq v^{\prime}$. Assume therefore that $x_{2}=v^{\prime}$,
otherwise $\{u^{\prime},v^{\prime}\}\cap Z=\varnothing$ and we are done. Now
$A=\{x_{1}^{\prime},x_{2},x_{3},x_{4}^{\prime},x_{5},...,x_{r}\}$ and
$y^{\prime}\nsim x_{2}=v^{\prime}$. If $x_{1}\in\operatorname{EPN}%
(x_{3},A^{\prime})$, then, by Theorem \ref{Thm_flip-set} applied to the
flip-set of $A^{\prime}$ using $\{x_{1},x_{1}^{\prime},x_{2},x_{4}^{\prime}%
\}$, this set is an $\operatorname{IR}(G)$-set that differs from $A^{\prime}$
in four elements. But this is impossible because $\operatorname{diam}(T)=3$.
Hence $x_{1}\notin\operatorname{EPN}(x_{3},A^{\prime})$ and the only
possibility is that $x_{1}\sim y^{\prime}$. Moreover, by Lemma
\ref{Lem_Not_Path}$(i)$ applied to $A^{\prime}$ and $A$, and because
$x_{1}^{\prime}\nsim x_{2}$, we see that $x_{1}^{\prime}\sim x_{4}^{\prime
}\sim x_{2}$. Consider the set $Q_{5}=\{x_{1},y^{\prime},x_{3},...,x_{r}\}$
and note that $G[Q_{5}]$ contains the path $(x_{3},x_{1},y^{\prime},x_{4})$.
Hence, by Lemma \ref{Lem_Number_IR-sets1}, $Q_{5}$ is not an
$\operatorname{IR}(G)$-set. But $x_{2}\in\operatorname{EPN}(x_{1},Q_{5})$,
$x_{2}^{\prime}\in\operatorname{EPN}(y^{\prime},Q_{5})$, $x_{3}^{\prime}%
\in\operatorname{EPN}(x_{3},Q_{5})$ and, unless $x_{1}\sim x_{4}^{\prime}$, we
also have that $x_{4}^{\prime}\in\operatorname{EPN}(x_{4},Q_{5})$; all
vertices $x_{i}$ with $i\geq5$ are isolated in $G[Q_{5}]$. Therefore
$x_{1}\sim x_{4}^{\prime}$. Finally, consider the set $Q_{6}=\{x_{4}^{\prime
},x_{2},...,x_{r}\}$, and note that $x_{1}^{\prime}\in\operatorname{EPN}%
(x_{4}^{\prime},Q_{6})$, $x_{2}^{\prime}\in\operatorname{EPN}(x_{2},Q_{6})$
and $y^{\prime}\in\operatorname{EPN}(x_{4},Q_{6})$, while $x_{i}$, for $i=3$
or $i\geq5$, is isolated in $G[Q_{6}]$. Hence $Q_{6}$ is an $\operatorname{IR}%
(G)$-set. However, $Q_{6}\overset{x_{4}^{\prime}x_{1}}{\sim}_{T}X$, from which
it follows that $T$ either has diameter at least $4$ (if $Q_{6}$ is
nonadjacent to all other $\operatorname{IR}(G)$-sets), or a cycle (otherwise).
This concludes the proof of Case~2.2 ($u=x_{j}$, $j\geq4$) and also of Case
2.~$\lozenge$\smallskip

\noindent\textbf{Case 3:\hspace{0.1in}}$y=x_{3}$, i.e., $A^{\prime}%
=\{x_{1}^{\prime},x_{2}^{\prime},y^{\prime},x_{4},\dots,x_{r}\}$ and $x_{3}%
\in\operatorname{EPN}(y^{\prime},A^{\prime})$. Similar to Case 2, if
$y^{\prime}\sim x_{3}^{\prime}$, then $X_{2}\overset{x_{3}y^{\prime}}{\sim
}A^{\prime}\overset{y^{\prime}x_{3}^{\prime}}{\sim}X^{\prime}\overset{x_{3}%
^{\prime}x_{3}}{\sim}X_{2}$, thus forming a cycle. Hence $y^{\prime}\nsim
x_{3}^{\prime}$. Recall that if the external $A^{\prime}$-private neighbours
$u^{\prime}$ and $v^{\prime}$ of $u$ and $v$, respectively, satisfy
$\{u^{\prime},v^{\prime}\}\cap Z=\varnothing$, then $d_{T}(X_{i},A)\geq2$ for
$i=1,2$, contrary to Observation \ref{Diam3_tree}. The possibilities for
$u^{\prime}$ and $v^{\prime}$ among $x_{i},x_{i}^{\prime}$, $i=1,2,3$, are
$x_{1},x_{2}$ and $x_{3}^{\prime}$. However, $x_{3}^{\prime}$ is adjacent to
both $x_{1}^{\prime}$ and $x_{2}^{\prime}$, so we only need to consider
$x_{1}$ and $x_{2}$. Note that
\begin{equation}
\text{if\ (say)\ }u^{\prime}=x_{1}\text{,\ then\ }u=x_{1}^{\prime
}\text{\ since\ }x_{1}\sim x_{1}^{\prime}\text{;\ similarly,\ if\ }v^{\prime
}=x_{2}\text{,\ then\ }v=x_{2}^{\prime}. \label{eq_case3a}%
\end{equation}

Suppose $u^{\prime}=x_{1}$ and $v^{\prime}=x_{2}$. Then $u=x_{1}^{\prime}$ and
$v=x_{2}^{\prime}$, hence, by Claim \ref{Cl2} and the elements in $X$ and
$A^{\prime}$, $x_{3}$ is the only neighbour of $y^{\prime}$ in $X$, i.e.,
$\{y^{\prime},x_{3}^{\prime}\}\subseteq\operatorname{EPN}(x_{3},X)$. But then,
since $y^{\prime}\nsim x_{3}^{\prime}$, $(X-\{x_{2},x_{3}\})\cup\{y^{\prime
},x_{3}^{\prime}\}$ is an independent set of $G$ of cardinality
$\operatorname{IR}(G)$, hence an independent $\operatorname{IR}(G)$-set
different from both $X_{1}$ and $X_{2}$, contrary to (\ref{eq_indep}).

Suppose $u=x_{1}^{\prime}$ and $v^{\prime}\neq x_{2}$, that is, $u^{\prime}$
is either $x_{1}$ or $u^{\prime}\notin Z$, and $v^{\prime}\notin Z$. Then
$A=\{u^{\prime},v^{\prime},x_{3},...,x_{r}\}$, so that, by the conditions on
$u^{\prime}$, the set $\{u^{\prime},v^{\prime}\}\subseteq A-(X_{1}\cup X_{2}%
)$. Consequently, $d_{T}(X_{i},A)\geq2$ for $i=1,2$, contrary to
Observation~\ref{Diam3_tree}. Similarly, if $u\neq x_{1}^{\prime}$ (so
$u^{\prime}\notin Z$) and $v^{\prime}\notin Z$, then $d_{T}(X_{i},A)\geq2$ for
$i=1,2$, giving the same contradiction. We examine the remaining case.

Suppose $u\neq x_{1}^{\prime}$ and $v^{\prime}=x_{2}$. Then $v=x_{2}^{\prime}$
and $u,u^{\prime}\notin Z$. Assume without loss of generality that $u=x_{4}$.
Then $A=\{x_{1}^{\prime},x_{2},x_{3},u^{\prime},x_{5},\dots,x_{r}\}$. We
consider the edges of $G[A\cup Z]$ and $G[A^{\prime}\cup Z]$. Applying Lemma
\ref{Lem_Not_Path}$(i)$ to $A^{\prime}$ and $A$ (and noting that $x_{2}\nsim
x_{3}$ and $x_{2}\nsim x_{1}^{\prime}\nsim x_{3}$) we find that $x_{2}\sim
u^{\prime}\sim x_{3}$.\emph{ }Since $u^{\prime}\in\operatorname{EPN}%
(x_{4},A^{\prime})$, $u^{\prime}$ is nonadjacent to all vertices in
$A^{\prime}-\{x_{4}\}=\{x_{1}^{\prime},x_{2}^{\prime},y^{\prime}%
,x_{5},...,x_{r}\}$. Since each of $x_{1}^{\prime},x_{2}^{\prime}%
,x_{3}^{\prime}$ has an $X^{\prime}$-external private neighbour, Lemma
\ref{Lem_Number_IR-sets1} implies that $\operatorname{EPN}(x_{4},X^{\prime
})=\varnothing$. The only remaining possibility is that $u^{\prime}\sim
x_{3}^{\prime}$. To avoid $A^{\prime}$ having four vertices with external
private neighbours, $x_{1}\notin\operatorname{EPN}(x_{1}^{\prime},A^{\prime}%
)$, and since $x_{2}^{\prime}\nsim x_{1}\nsim x_{4}$, necessarily $y^{\prime
}\sim x_{1}$.

\begin{figure}[ptb]%
\centering
\includegraphics[
height=1.9424in,
width=4.9139in
]%
{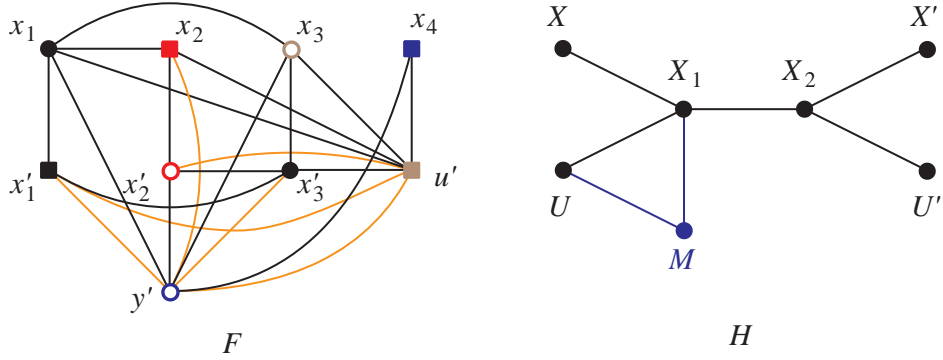}%
\caption{The graph $F=G[Z\cup\{x_{4},u^{\prime},y^{\prime}\}]$ and the graph
$H$ formed by some of its $\operatorname{IR}$-sets, as described in Case 3 of
the proof of Theorem \ref{Thm_Main}.}%
\label{Fig_Case2.3}%
\end{figure}

Consider the set $Q_{7}=\{x_{1},x_{2},x_{2}^{\prime},y^{\prime},x_{5}%
,...,x_{r}\}$. Since $x_{1}\sim x_{2}\sim x_{2}^{\prime}=v\sim y^{\prime}$,
$x_{1},x_{2},x_{2}^{\prime},y^{\prime}$ all have positive degree in $G[Q_{7}%
]$. Hence $Q_{7}$ is not an $\operatorname{IR}(G)$-set (by Lemma
\ref{Lem_Number_IR-sets1}). However, if $x_{1}\nsim u^{\prime}$, then
$x_{1}^{\prime}\in\operatorname{EPN}(x_{1},Q_{7}),\ u^{\prime}\in
\operatorname{EPN}(x_{2},Q_{7}),\ x_{3}^{\prime}\in\operatorname{EPN}%
(x_{2}^{\prime},Q_{7})$ and $x_{4}\in\operatorname{EPN}(y^{\prime},Q_{7})$.
Thus $x_{1}\sim u^{\prime}$. The subgraph $F=G[Z\cup\{x_{4},u^{\prime
},y^{\prime}\}]$ is shown by the black edges in Figure \ref{Fig_Case2.3},
where the edges incident with $u^{\prime}$ and $y^{\prime}$ in the complement
$\overline{F}$ are shown in colour (light grey in monotone). Recall that
$U=\{x_{1},x_{1}^{\prime},x_{2},x_{4},...,x_{r}\}$ and consider the set
$M=\{x_{1}^{\prime},x_{2},x_{4},u^{\prime},x_{4},...,x_{r}\}$. The vertices
$x_{2},x_{4}$ and $u^{\prime}$ have positive degree in $G[M]$ and have
$M$-external private neighbours $x_{2}^{\prime},y^{\prime}$ and $x_{3}$,
respectively, hence $M$ is an $\operatorname{IR}(G)$-set. However,
$M\overset{u^{\prime}x_{1}}{\sim}_{T}U\overset{x_{1}x_{3}}{\sim}_{T}%
X_{1}\overset{x_{3}u^{\prime}}{\sim}_{T}M$, hence $T$ contains a cycle, which
is impossible.%

We conclude that $\{u^{\prime},v^{\prime}\}\cap Z=\varnothing$, which leads to
a contradiction with Observation \ref{Diam3_tree} as explained
above.~$\lozenge$\smallskip

Therefore, no $y\in\{x_{1}^{\prime},x_{2}^{\prime},x_{3}\}$ is swapped when
obtaining $A^{\prime}$ from $X_{2}=\{x_{1}^{\prime},x_{2}^{\prime},x_{3}%
,\dots,x_{r}\}$; it follows that $y=x_{j}$ for some $j\geq4$.~$\blacklozenge
$\medskip

We have shown in Claims \ref{Cl3} and \ref{Cl4} that if $A^{\prime}%
\notin\{X_{1},X_{3},U^{\prime}\}$ is any $\operatorname{IR}(G)$-set such that
$X_{2}\overset{yy^{\prime}}{\sim}_{T}A^{\prime}$, then $y=x_{i}$ for $i\geq4$,
and $y^{\prime}\in V(G)-(X\cup Z)$. Assume without loss of generality that
$y=x_{4}$ and denote $y^{\prime}$ by $x_{4}^{\prime}$, so that $A^{\prime
}=\{x_{1}^{\prime},x_{2}^{\prime},x_{3},x_{4}^{\prime},x_{5}\dots,,x_{r}\}$
and $x_{4}\in\operatorname{EPN}(x_{4}^{\prime},A^{\prime})$. We determine the
two neighbours of $x_{4}^{\prime}$ in $A^{\prime}$ in Claim \ref{Cl5}.

\begin{claim}
\label{Cl5}The neighbours of $y^{\prime}=x_{4}^{\prime}$ in $A^{\prime}$ are
$x_{2}^{\prime}$ and a vertex $x_{j}\in X$, where $j\geq5$.
\end{claim}

\noindent\textbf{Proof of Claim \ref{Cl5}.\hspace{0.1in}}For the two
neighbours $u$ and $v$ of $x_{4}^{\prime}$ in $A^{\prime}$, we consider the
possibilities for their external $A^{\prime}$-private neighbours $u^{\prime}$
and $v^{\prime}$, respectively. Clearly, if $\{u^{\prime},v^{\prime
}\}\subseteq A-(X_{1}\cup X_{2})$, then $d_{T}(X_{i},A)\geq2$ for $i=1,2$ and
we obtain a contradiction to Observation \ref{Diam3_tree} as before. Therefore
we may assume that
\begin{equation}
|\{u^{\prime},v^{\prime}\}\cap A\cap X_{1}|\geq1\text{ \ or \ }|\{u^{\prime
},v^{\prime}\}\cap A\cap X_{2}|\geq1. \label{eq_u'v'}%
\end{equation}
By the properties of the $X_{i}$, no vertex in $\{x_{5},...,x_{r}\}$ is
adjacent to a vertex in $X\cup Z$. Thus, if (say) $u\in\{x_{5},...,x_{r}\}$,
then $u^{\prime}\notin X_{1}\cup X_{2}$, and the same holds for $v$ and
$v^{\prime}$. Hence the inequalities (\ref{eq_u'v'}) imply that $|\{u,v\}\cap
\{x_{5},...,x_{r}\}|\leq1$ and $|\{u,v\}\cap\{x_{1}^{\prime},x_{2}^{\prime
},x_{3}\}|\geq1$. Since $x_{1}^{\prime}\sim x_{1}\sim x_{3}$ and
$x_{1}^{\prime}\sim x_{3}^{\prime}\sim x_{3}$, we know that $\{u^{\prime
},v^{\prime}\}\cap\{x_{1},x_{3}^{\prime}\}=\varnothing$. Since we also know
that $u^{\prime},v^{\prime}\notin A^{\prime}$, the inequalities (\ref{eq_u'v'}%
) now reduce to $\{u^{\prime},v^{\prime}\}\cap\{x_{1},x_{2},x_{3}^{\prime
}\}=\{x_{2}\}$. Necessarily, then, $x_{4}^{\prime}$ is adjacent to
$x_{2}^{\prime}$ and to some vertex $x_{j}\in X$, where $j\geq5$, as
required.~$\blacklozenge$\medskip\ 

Assume without loss of generality that $u^{\prime}=x_{2}$, that is,
$u=x_{2}^{\prime}$, and that $v$ and $v^{\prime}$ are $x_{5}$ and
$x_{5}^{\prime}$, respectively. Then $A^{\prime}=\{x_{1}^{\prime}%
,x_{2}^{\prime},x_{3},x_{4}^{\prime},x_{5},...,x_{r}\}$ and $A=\{x_{1}%
^{\prime},x_{2},x_{3},x_{4},x_{5}^{\prime},...,x_{r}\}$. Let $B=(A-\{x_{4}%
\})\cup\{x_{5}\}$ and $B^{\prime}=(A^{\prime}-\{x_{5}\})\cup\{x_{4}\}$ be the
skip-sets of $A$ and $A^{\prime}$, respectively. By Lemma \ref{Lem_Not_Path},
$\mathcal{A}=\{A,A^{\prime},B,B^{\prime}\}$ is a $4$-cluster. Clearly,
$\mathcal{A}$ and $\mathcal{X}=\{X,X^{\prime},U,U^{\prime}\}$ are disjoint $4$-clusters.

We complete the proof of the theorem by showing, in Claim \ref{Cl6}, that any
other leaf of $T$ belongs to a $4$-cluster disjoint from both $\mathcal{A}$
and $\mathcal{X}$.

\begin{claim}
\label{Cl6}Any $\operatorname{IR}(G)$ set $C^{\prime}\notin\{X_{1},X^{\prime
},A^{\prime},B^{\prime},U^{\prime}\}$ adjacent to $X_{2}$ is obtained by
swapping a vertex $x_{j}\in X_{2}$ for a vertex $x_{j}^{\prime}$, where
$j\geq6$, $x_{j}^{\prime}\notin\{x_{1},...,x_{5},x_{1}^{\prime},...,x_{5}%
^{\prime}\}$, the only neighbours of $x_{j}^{\prime}$ in $C^{\prime}$ are
$x_{2}^{\prime}$ and $x_{k}$ for some$\ k\geq6,\ k\neq j$, and the $C^{\prime
}$-external private neighbour $x_{6}^{\prime}$ of $x_{6}$ also does not belong
to $\{x_{1},...,x_{5},x_{1}^{\prime},...,x_{5}^{\prime}\}$.
\end{claim}

\noindent\textbf{Proof of Claim \ref{Cl6}.\hspace{0.1in}}Suppose that
$C^{\prime}\notin\{X_{1},X^{\prime},A^{\prime},B^{\prime},U^{\prime}\}$ and
$X_{2}\overset{yy^{\prime}}{\sim}_{T}C^{\prime}$. We first consider
$C^{\prime}$ in relation to the $4$-cluster $\mathcal{X}$, repeating the
arguments above.

By Claim \ref{Cl2}, $y^{\prime}$ is adjacent to exactly two vertices of
$C^{\prime}$, and by Claims \ref{Cl3} and \ref{Cl4}, $y=x_{j}$, where $j\geq4$
and $y^{\prime}\in V(G)-(X\cup\{x_{1}^{\prime},x_{2}^{\prime},x_{3}^{\prime
}\})$. Say $y^{\prime}=x_{j}^{\prime}$. By Claim \ref{Cl5}, the neighbours of
$x_{j}^{\prime}$ in $C^{\prime}$ are $x_{2}^{\prime}$ and a vertex $x_{k}\in
X$, where $k\geq4$ and $k\neq j$. \label{hereFri}Note that $G[\{x_{2}%
,x_{2}^{\prime},x_{4},x_{4}^{\prime},x_{5},x_{5}^{\prime}\}]\cong G[Z]$ under
an isomorphism that fixes $x_{2}$ and $x_{2}^{\prime}$. Repeating the proof
above for the $4$-cluster $\mathcal{A}$, Claims \ref{Cl3} and \ref{Cl4} give
that $x_{j}\in\{x_{1},x_{1}^{\prime},x_{3},x_{3}^{\prime},x_{6},...,x_{r}\}$
and $x_{j}^{\prime}\in V(G)-(X\cup\{x_{2}^{\prime},x_{4}^{\prime}%
,x_{5}^{\prime}\})$. The two conditions on $x_{j}$ show that $j\geq6$, while
those on $x_{j}^{\prime}$ show that $x_{j}^{\prime}\in V(G)-(X\cup
\{x_{1}^{\prime},...,x_{5}^{\prime}\})$. Hence assume $j=6$, i.e., $y=x_{6}$
and $y^{\prime}=x_{6}^{\prime}$. Now Claim \ref{Cl5} applied to $\{x_{2}%
,x_{2}^{\prime},x_{4},x_{4}^{\prime},x_{5},x_{5}^{\prime},x_{6},x_{6}^{\prime
}\}$ asserts that the neighbours of $x_{6}^{\prime}$ in $C^{\prime}$ are
$x_{2}^{\prime}$ and a vertex in $C^{\prime}$ different from $x_{2}^{\prime
},x_{4},x_{5}$. Consolidating the two conditions obtained from Claim
\ref{Cl5}, the neighbours of $x_{6}^{\prime}$ in $C^{\prime}$ are
$x_{2}^{\prime}$ and $x_{k}$, where $k\geq7$; say $k=7$. Let $x_{7}^{\prime
}\in\operatorname{EPN}(x_{7},C^{\prime})$. As above, $G[\{x_{2},x_{2}^{\prime
},x_{6},x_{6}^{\prime},x_{7},x_{7}^{\prime}\}]\cong G[Z]$ under an isomorphism
that fixes $x_{2}$ and $x_{2}^{\prime}$. Thus $(x_{6},x_{6}^{\prime}%
,x_{7},x_{7}^{\prime},x_{6})$ is a $4$-cycle in $G$ in which $x_{6}$ and
$x_{7}$ are nonadjacent vertices of degree $2$, and the skip-set $D^{\prime
}=(C^{\prime}-\{x_{7}\})\cup\{x_{6}\}$ and the flip-sets $C$ and $D$ of
$C^{\prime}$ and $D^{\prime}$ obtained by using the obvious private neighbours
are all $\operatorname{IR}(G)$-sets which form a $4$-cluster.~$\blacklozenge
$\medskip

By repeating the arguments in the proof of Claim \ref{Cl6}, considering all
existing $4$-clusters in turn in each step until all vertices in $X$ have been
used, we show that the leaves of $T$ occur in disjoint $4$-clusters, thus
obtaining the desired result.~$\blacksquare$

\section{Open Problems}

\label{Sec_Open}We conclude by mentioning a number of conjectures and open
problems, most of which also appear in \cite{mynhardt19}.

\begin{conjecture}
\emph{\cite{mynhardt19}}$\hspace{0.1in}(i)\hspace{0.1in}$The path $P_{n}$ is
not an $\operatorname{IR}$-graph for each $n\geq3$.

\begin{enumerate}
\item[$(ii)$] The cycle $C_{n}$ is not an $\operatorname{IR}$-graph for each
$n\geq5$.
\end{enumerate}
\end{conjecture}

\begin{problem}
\emph{\cite{mynhardt19}\hspace{0.1in}}Prove or disprove: Complete graphs and
$K_{m}\boksie K_{n}$, where $m,n\geq2$, are the only connected claw-free
$\operatorname{IR}$-graphs.
\end{problem}

\begin{problem}
\emph{\cite{mynhardt19}\hspace{0.1in}}Determine which double spiders are
$\operatorname{IR}$-trees.
\end{problem}

\begin{problem}
\emph{\cite{mynhardt19}\hspace{0.1in}}Characterize $\operatorname{IR}$-graphs
having diameter $2$.
\end{problem}

\begin{problem}
Determine which graphs are $\operatorname{IR}$-graphs of a graph $G$ such that
$\operatorname{IR}(G)=2$.
\end{problem}

\noindent\textbf{Acknowledgements\hspace{0.1in}}We acknowledge the support of
the Natural Sciences and Engineering Research Council of Canada (NSERC), PIN 253271.

\noindent Cette recherche a \'{e}t\'{e} financ\'{e}e par le Conseil de
recherches en sciences naturelles et en g\'{e}nie du Canada (CRSNG), PIN
253271.%
\begin{center}
\includegraphics[
natheight=0.773100in,
natwidth=1.599900in,
height=0.4151in,
width=0.8276in
]%
{../NSERC2014_2019/NSERC_DIGITAL_BW/NSERC_DIGITAL_BW/NSERC_BLACK.jpg}%
\end{center}

\end{document}